# THE ROUSSEAU NUMBER: AN INFORMETRIC VERSION OF THE ERDÖS NUMBER

## LEO EGGHE

## HASSELT UNIVERSITY, BELGIUM

Leo.egghe@uhasselt.be

ORCID : 0000-0001-8419-2932

## ABSTRACT

Some mathematical properties of the Erdös number and its formal equivalents are shown.  An informetric equivalent is presented: the Rousseau number. This contribution honors Ronald Rousseau for his 73$^{th}$birthday.

## INTRODUCTION

Paul Erdös (1913-1996) was a brilliant mathematician who collaborated extensively with colleagues, resulting in more than 1,500 papers and more than 500 different co-authors. This fact was the basis of the idea for the so-called Erdös number, a natural number associated with a person and defined as follows (see e.g. Glänzel and Rousseau (2005)) : Erdös himself has Erdös number 0. All of his co-authors have Erdös number 1. All authors who co-authored a publication with an author with Erdös number 1 receive Erdös number 2, provided they do not have the Erdös number 1, and so on.

For example my Erdös number is 3 because I co-authored with Bellow who co-authored with Dvoretzky who co-authored with Erdös and no shorter paths exist. Erdös numbers are not limited to mathematicians. For instance, Albert Einstein has Erdös number 2. Lists of mathematicians, physicists, chemists and even medicine scientists and their Erdös numbers can be found on the Internet, see https://en.wikipedia.org/wiki/List_of_people_by_Erd%C5%91s_number

It is somewhat surprising that there exists an equivalent of the Erdös number in the field of motion pictures: the Bacon number, relating to the American actor Kevin Bacon. The same definition as the Erdös number can be given for the Bacon number, now replacing co-authors by co-actors. One can make a critical remark here: unlike mathematics (where only co-authors make a mathematical publication), in motion pictures, a film has more collaborators than just actors (e.g. the director,…), and hence, instead of only using actors in the definition of the Bacon number, one should use (e.g.) the names in the title of a motion picture for a more natural definition of the Bacon number. Since this is outside the scope of this paper (and journal) we leave it here as it is.

It is well-known that all these numbers are "surprisingly" small (i.e., smaller than one would expect). This is linked with the so-called "Small World Phenomenon" in social networks where (e.g., in the world) "strangers are linked through a short chain of acquaintances" hereby finding that, for any two persons A and B in the world one needs, on the average, only chains of acquaintances of length 6 to make the connection between A and B, called "Six Degrees of Separation" – see also Kleinberg (2004), Rousseau, Egghe and Guns (2018).

Such A-numbers (relating to a person A) can, of course, be defined for any person A, provided A has links (as e.g. described above) to other persons and this can be defined in any field, besides the already mentioned fields of mathematics and motion pictures. Since A-numbers, in general, can be defined in networks (graphs), as studied e.g. in informetrics, it is therefore surprising that there has not been defined an A-number for an informetrician A. This will be done in the last section of this note, but first, in the next section, we will define the A-number in general (weighted) networks and we will also present some elementary properties.

**THE A-NUMBER IN NETWORKS**

A network (see Egghe and Rousseau (1990), Egghe (2005), Rousseau, Egghe and Guns (2018) – just to mention 3 books in informetrics) – also called a graph – is a pair of finite sets (V,E) consisting of points (vertices, nodes) in V and links (edges, arcs) between vertices, say (A,B) $\in$ E (in this order) with A,B $\in$ V. In this note we will mainly deal with social networks, where e.g. nodes are persons and edges are links between them (friends, co-author (or actor) – ship, …). But nodes in a social network can also be papers (linked e.g., through references or citations) or even cities linked by roads as edges.

We will say that a network is directed if, for A,B $\in$ V, the order in the edge (A,B) $\in$ E matters. Example : A cites B (in a paper) and then we have the edge (A,B) in this order (since it differs from (B,A) which might not even exist as e.g. in the case that A and B are articles). If the order in the edge (A,B) does not matter (mathematically : (A,B) = (B,A)) we talk about an undirected network. Example : A is co-author of B, hence B is also co-author of A. We can also talk of a collaboration network. Henceforth in this note, we will only work with this type of network (V,E).

If, to each node (A,B) $\in$ E, a value w(A,B) > 0 is attached (with w : E $\rightarrow$ $R^+$), we say that the network (V,E) is weighted. These weights can be 1 in which case we talk about an unweighted network (as a special example of a weighted network !) : therefore, all our assertions for weighted networks also apply to unweighted networks. A path in a network from A $\in$ V to B $\in$ V is a sequence in V : A=C(0), C(1), … C(n-1), C(n)=B and edges (C(i-1),C(i)) $\in$ E, $\forall$i=1,…,n with n $\in$ N. For n=1 we have that (A,B) $\in$ E. A network is called connected if, for all A,B $\in$ V, there exists at least one path from A to B. We henceforth suppose that we have a connected network (and will omit the adjective "connected"). This is an evident assumption in collaboration studies one can only look at connected networks (or at the connected subgraphs in case the graph (= network) is not connected). Denote w(i)= w(C(i-1),C(i)), the weight of the edge (C(i-1),C(i)) in the above path. Then we say that $\sum_{i=1}^{n} w(i)$ is the length of this path. In an unweighted network this length is simply n = the number of edges from A to B in this path. We define a shortest path between A and B as a path of minimum length. A synonym is a network geodesic between A and B. Obviously such a network geodesic is not necessarily unique but its length is, called the geodesic distance between A and B, denoted d(A,B) and where we define d(A,A)=0, for all A ∈ V. That d is indeed a distance on V is proved in the next Proposition (see also Goddard and Oellermann (2010)).

**Proposition 1**: d is an metric (distance) and hence V,d, is a metric space.

Proof : We have to prove the following three properties of a metric (for A,B,C ∈ V) :

  (i)      d(A,B) = 0 if and only if A = B
  (ii)     d(A,B) = d(B,A)
  (iii)    d(A,B) ≤ d(A,C) + d(C,B), the so-called triangle inequality.

For (i), if A=B, then d(A,B) = d(A,A) = 0 by definition. If A ≠ B then any path between A and B (which exists by assumption) has at least one edge with a positive weight, hence d(A,B) > 0 being at least the minimum of these finite number of positive weights. (ii) is clear since (A,B) = (B,A) in a collaboration network and since w is a function on E. For (iii) we have that, going from A to B via C is one way of making a path from A to B. Since d(A,B) is the length of the smallest path from A to B, (iii) is clear. □

**Definition** : In a network (V,E), for A ∈ V, we define the A-number of B ∈ V as the geodesic distance d(A,B) between A and B.

Note that this definition depends on the weights in the network. So, here, we give an explicit definition of d(A,B) :

d(A,B) = min { $\sum_{i=1}^{n} w(i)$ ∥ (C(j))$_{j=0,\dots,n}$ is a path between A and B in V }          (1)

(hence C(0) = A and C(n) = B).

We note the following trivial Lemma :

**Lemma** : for all A,B ∈V :

d(A,B) = min { $\sum_{i=1}^{n} d(C(i-1), C(i))$ ∥ (C(j))$_{j=0,\dots,n}$ is a path between A and B }          (2)

Proof : ≥ in (2) is clear since, for all i=1,…,n : d(C(i-1),C(i)) ≤ w(i) since the path (C(i-1),C(i)) with length w(i) is only one path between C(i-1) and C(i). Now we prove the opposite inequality ≤ in (2). For every path C(0)=A, C(1), …, C(n-1), C(n)=B between A and B and for all i=1,…,n, we – possibly – extend this path by adding vertices between C(i-1) and C(i) that yield a geodesic between C(i-1) and C(i). Denoting this new path again by $(C(j))_{j=0,\dots n}$ we have that its length is given by $\sum_{i=1}^{n} d(C(i-1), C(i))$. Hence, by definition of d(A,B) : ≤ is proved in (2). Hence (2) is proved. □

**Examples**

(i) If we take w(i) = 1 for all i=1,…n we obtain

$A_B$=: d(A,B) = min {n ‖ $(C(j))_{j=0,\dots n}$ is a path between A and B }　　　(3)

being the classical definition of e.g. the Erdös number of B, in case A = Erdös =: E

(ii) Now we replace w(i) = 1 in (i) by

$$w(i) = \frac{1}{\#(C(i-1),C(i))}　　　(4)$$

for i=1,…n : 1 divided by the number of links in the edge (C(i-1),C(i)). For collaboration this means the number of times C(i-1) and C(i) have collaborated with each other (e.g. co-authored publications). This is a very logical weighted variant of (i) since it takes the strength of collaboration into account. Now, by (1) and (4) :

$A'_B$ =: d'(A,B) = min { $\sum_{i=1}^{n} \frac{1}{\#(C(i-1),C(i))}$ ‖ $(C(j))_{j=0,\dots n}$ is a path between A and B }　　　(5)

We have the following easy Proposition :

**Proposition 2** : For all A,B,C ϵV :

(i)

$A_B$ = 0 ⟺ A = B　　(6)

$A_B$ = $B_A$　　　(7)

$A_B ≤ A_C + C_B$　　(8)

$A_B ≤ A_C + B_C$　　(9)

(ii)

$A'_B$ = 0 ⟺ A = B　　(6')

$A'_B$ = $B'_A$　　　(7')

$A'_B ≤ A'_C + C'_B$　　(8')

$A'_B ≤ A'_C + B'_C$　　(9')

(iii)

$A'_B ≤ A_B$　　　(10)

$$A'_B \leq \frac{1}{\#(A,B)} \qquad (11)$$

Proof : (i) and (ii) follow from Proposition 1. (10) follows since in (4), $\#(C(i-1),C(i)) \geq 1$. Finally (11) follows by the following argument : Suppose first that A and B did not collaborate directly. Then $\#(A,B) = 0$ and (11) is trivial. Suppose now that A and B collaborated directly. Then $\frac{1}{\#(A,B)}$ is one number in (5) of which $A'_B$ is the minimum (since (A,B) is one path from A to B). Note that in this particular case, $A'_B \leq \frac{1}{\#(A,B)} \leq 1$.

□

To make a distinction between $A_B$ and $A'_B$ in words one could say that $A_B$ is the classical, unweighted A-number of B while $A'_B$ is the weighted A-number of B.

Illustration : A has 2 publications co-authored with B and B has 4 publications co-authored with C. Then $A_B = B_C = 1$, $A_C = 2$, $A'_B = \frac{1}{2}$, $B'_C = \frac{1}{4}$ and $A'_C = \frac{1}{2} + \frac{1}{4}$ in which case (8), (9), (8') and (9') are equalities.

Remark : (10) and (11) cannot be combined for two reasons : $\frac{1}{\#(A,B)} \leq A_B$ is false if $\#(A,B) = 0$ and $A_B$ exists (e.g. for A = Erdös, B = Egghe). But also $A_B \leq \frac{1}{\#(A,B)}$ is false in any case where $\#(A,B) > 1$ (since $A_B = 1$).

**The Rousseau number**

Ronald Rousseau is a well-known prolific informetrician and mathematician who has published just over 500 official publications, a top number in informetrics. What is even more top is his number of different co-authors in these publications which is larger than 250. Also, around 350 of his publications are co-authored (i.e. around 70 %). These high numbers were the basis of the idea of your author to – officially – introduce the Rousseau number in our field (officially since, in principle, anyone can be attributed an own number in a similar way as the Erdös number is attributed to Paul Erdös). So we hereby define the Rousseau number $R_A$ of author A as $d(R,A) = R_A$ as in (3) and the weighted Rousseau number $R'_A$ of an author A as $d'(R,A) = R'_A$ as in (5).

So, $R_R = 0$ and $R_A = 1$ for all authors A who co-authored with Rousseau, including your author. Reserving E for Erdös and Eg for Egghe we hence have $R_{Eg} = 1$ but, since Egghe and Rousseau collaborated so much (resulting in (at this time) 73 joint publications) we have that $R'_{Eg} = 0.014$ bringing Egghe very close to Rousseau without ever being at a distance 0 (this is only the case for Rousseau himself).

Using formulae (8), (9), (8') and (9'), we can establish relationships between the Rousseau number and (e.g.) the famous Erdös number. Since $E_{Eg} = 3$ and $R_{Eg} = 1$ we have that $E_R \leq E_{Eg} + R_{Eg} = 3 + 1 = 4$ (by (9)) so that Rousseau's Erdös number is (at most) 4 (and probably equal to 4). Hence all the authors A who collaborated with Rousseau directly have Erdös number $E_A \leq 5$. We must acknowledge that the calculation of the exact $A'_B$ numbers is much more complicated than the $A_B$ ones since we need to know, for all collaborators in every edge between A and B, the number of joint publications ! But even if we do not have the complete data, the inequalities (8') and (9') and (10) are handy for an upper estimate of $A'_B$. Indeed, since Egghe has 2 publications with Bellow we have that, by (5), $E'_{Eg} \leq 2.5$ only using that Bellow and Dvoretzky and Dvoretzky and Erdös have at least one publication together. It then follows from (9') and (10) that $E'_R \leq E'_{Eg} + R_{Eg} \leq 3.5$ without more work. Of course, using the extra knowledge that $R'_{Eg} = 0.014$ we obtain, by (9') that $E'_R \leq E'_{Eg} + R'_{Eg} \leq 2.514$. Hence, without further work, we have that $E'_A \leq 3.514$ for all co-authors of Rousseau.

Your author considers the Introduction of the Rousseau number in informetrics as an original and well-deserved present for Rousseau's 73[th]birthday. Of course, if the reader prefers 75, he/she might remember this paper in 2024. I preferred 73 above 75 since 73 is the closest prime number to 75 and since, for mathematicians, prime numbers are superior to composite numbers such as 75.